\documentclass[11pt,a4paper,leqno]{amsart}
\usepackage{amssymb, amsmath, amscd}
\usepackage[all]{xy}
\usepackage{color}

\newtheorem{theorem}{Theorem}[section]

\newtheorem{proposition}[theorem]{Proposition}
\newtheorem{remark}[theorem]{Remark}

\newcommand{\bach}{\mathfrak{B}}
\newcommand{\Id}{\operatorname{Id}}
\newcommand{\Ric}{\operatorname{Ric}}
\newcommand{\Hes}{\operatorname{Hes}}

\begin{document}
\title[Conformal geometry of non-reductive homogeneous spaces]
{Conformal geometry of non-reductive four-dimensional homogeneous spaces}
\author[Calvi\~{n}o-Louzao, Garc\'ia-R\'io, Guti\'errez-Rodr\'\i guez,  V\'{a}zquez-Lorenzo]{E. Calvi\~{n}o-Louzao, E. Garc\'ia-R\'io, I. Guti\'errez-Rodr\'\i guez,  \\
R. V\'{a}zquez-Lorenzo}
\address{(E. C.-L.) Department of Mathematics, IES Ram\'{o}n Caama\~{n}o, Mux\'ia, Spain}
\address{(E. G.-R.) Faculty of Mathematics,
University of Santiago de Compostela, 15782 Santiago de Compostela, Spain}
\address{(I. G.-R.)  Faculty of Mathematics,
University of Santiago de Compostela, 15782 Santiago de Compostela, Spain}

\address{(R. V.-L.) Department of Mathematics, IES de Ribadeo Dionisio Gamallo, Ribadeo, Spain}

\email{estebcl@edu.xunta.es, eduardo.garcia.rio@usc.es, dzohararte@hotmail.fr, 
ravazlor@edu.xunta.es}
\thanks{Supported by projects GRC2013-045, MTM2013-41335-P and EM2014/009 with FEDER funds (Spain).}
\date{}

\begin{abstract}
We classify non-reductive  four-dimensional homogeneous conformally Einstein manifolds.

\end{abstract}

\maketitle

\section{Introduction}

Einstein metrics have a privileged role in pseudo-Riemannian geometry. 
Several generalizations of the Einstein condition have been intensively studied in the literature, Ricci solitons and quasi-Einstein metrics being just some examples.
The fact that the Ricci tensor is not preserved by conformal transformations motivated the study of conformally Einstein metrics, as those pseudo-Riemannian manifolds $(M,g)$ which contain an Einstein representative in the conformal class $[g]$. 

Brinkmann \cite{Brinkmann24} determined the necessary and sufficient conditions for a manifold to be conformally Einstein, which are codified in the conformally Einstein equation
$$
(n-2)\Hes_\varphi +\varphi\,\rho= \theta g\,,
$$
for some function $\theta$ on $M$. 
Besides its apparent simplicity, the integration of the conformally Einstein equation is surprisingly difficult. The equation is trivial in dimension two and it is equivalent to local conformal flatness in dimension three, which shows that it is overdetermined in the generic situation. Hence, the first non-trivial case to study is dimension four, where harmonicity of the Weyl tensor is a necessary condition to be conformally Einstein. 

Gover and Nurowski \cite{GoNu} obtained some tensorial obstructions
for a metric to be conformally Einstein under some non-degeneracy conditions for the conformal Weyl tensor. 
It is still an open question to characterize conformally Einstein manifolds by tensorial equations. In dimension four the Bach flat equation is only a necessary condition since there are Bach flat spaces which are not conformally Einstein. An important class of four-dimensional conformally Einstein metrics is obtained by considering the product of surfaces with nowhere zero scalar curvature under some additional conditions (see the survey \cite{Kuhnel-Rademacher, KR2} for more information). Moreover in such cases the conformal Einstein metric is unique up to a constant. 

The purpose of this paper is to analyze the conformally Einstein equation for a class of strictly pseudo-Riemannian four-dimensional homogeneous spaces, namely the non-reductive ones. We determine explicitly which non-reductive homogeneous four-manifolds are conformally Einstein and give all the possible conformal Einstein metrics in each case. 
It is worth to remark that all Einstein metrics inside each conformal class are Ricci flat
and,  moreover, they are not unique depending on the cases, allowing the existence of two-parameter and three-parameter families of conformal Ricci flat metrics in some cases.

The paper is organized as follows. The classification of the non-reductive four-dimensional homogeneous spaces given in \cite{Fels-Renner06} and the local forms of the metrics corresponding to the different classes obtained in \cite{Calvaruso-Fino-Zaeim15} are briefly reviewed in Section~\ref{se:2-1}. The Bach tensor is introduced in Section~\ref{se:2-2} and the classification of all Bach flat non-reductive four-dimensional homogeneous spaces is given in Theorem~\ref{th: Bach flat}. The conformally Einstein equation is treated in Section~\ref{se:2-3} where  the Main Theorem is stated, classifying the conformally Einstein non-reductive four-dimensional homogeneous spaces.
All the curvature calculations are carried out in Section~\ref{se:3}, while the proof of the Main Theorem is given in Section~\ref{se:4}.

\section{Preliminaries}\label{se:preliminaries}

Let $(M,g)$ be a connected pseudo-Riemannian manifold, $\nabla$ its   Levi-Civita connection and $R(X,Y)=\nabla_{[X,Y]}-[\nabla_X,\nabla_Y]$ the corresponding curvature tensor.
The Ricci tensor and the Ricci operator are given by $\rho(X,Y)=~g(\Ric\,X,Y)$ $=$ $\operatorname{trace}\{Z \mapsto R(X,Z)Y\}$. 
The scalar curvature $\tau=\operatorname{trace}\Ric$ is the metric trace of the Ricci tensor.


The curvature of any pseudo-Riemannian manifold $(M,g)$, as an endomorphism of the space of 2-forms $\Lambda^2=\Lambda^2(M)$, naturally decomposes under the action of the orthogonal group as
$R=\frac{\tau}{n(n-1)}\operatorname{Id}_{\Lambda^2} + \rho_0+W$, where $n=\dim M$, $\rho_0=\rho-\frac{\tau}{n}g$  is the trace-free Ricci tensor and $W$ is the Weyl curvature tensor.

\subsection{Non-reductive homogeneous spaces}\label{se:2-1}
A pseudo-Riemannian manifold is homogeneous if there is a group of isometries which acts transitively on $M$. Let $G$ be such a group of isometries and let $H$ denote the isotropy group at some fixed point. Then $(M,g)$ can be identified with the quotient space $(G/H,\tilde g)$, where $\tilde g$ is an invariant metric on $G$.
A homogeneous space $G/H$ is said to be \emph{reductive} if the Lie algebra admits a decomposition of the form $\mathfrak{g}=\mathfrak{h}\oplus\mathfrak{m}$ where $\mathfrak{m}$ is an $Ad(H)$-invariant complement of $\mathfrak{h}$. While every Riemannian homogeneous space is reductive, there are pseudo-Riemannian homogeneous spaces without any reductive decomposition. The geometry of reductive pseudo-Riemannian manifolds presents some similarities with the Riemannian case (see, for example \cite{GaOu}), but little is known about the non-reductive case.
The geometry of non-reductive homogeneous spaces is therefore an important aspect towards a good understanding of pseudo-Riemannian homogeneous manifolds.

Any homogeneous pseudo-Riemannian manifold is reductive in dimension two and three. In dimension four, a complete classification of non-reductive homogeneous spaces was obtained in  \cite{Fels-Renner06}. Later on a coordinate description was given in \cite{Calvaruso-Fino-Zaeim15} which we recall in order to state our results.

\begin{theorem} 
\label{th:1-1}
Let $(M,g)$ be a non-reductive homogeneous pseudo-Riemann\-ian manifold of dimension four. Then it is locally isometric to one of the following:
\begin{enumerate}
\item[(A.1)] $\mathbb{R}^4$ with coordinates $(x_1,x_2,x_3,x_4)$ and metric tensor

\smallskip\noindent
$\displaystyle
\begin{array}{l}
g=(4b x_2^2+a)\,dx_1^2+4bx_2\,dx_1dx_2-(4ax_2x_4-4cx_2+a)\,dx_1dx_3 \\
\noalign{\medskip}
\phantom{g}+4ax_2\,dx_1dx_4+b\,dx_2^2-2(ax_4-c)\,dx_2dx_3+2a\,dx_2dx_4+q\,dx_3^2\,,
\end{array}
$

\smallskip\noindent
where $a$, $b$, $c$ and $q$ are arbitrary constants with $a(a-4q)\neq 0$.

\medskip
\item[(A.2)] $\mathbb{R}^4$ with coordinates $(x_1,x_2,x_3,x_4)$ and metric tensor

\smallskip\noindent
$\displaystyle
\begin{array}{l}
g =  -2 a e^{2\alpha x_4}\,dx_1 dx_3+a e^{2\alpha x_4}dx_2^2+b e^{2(\alpha-1)x_4}dx_3^2 \\
\noalign{\medskip}
\phantom{g}+2c e^{(\alpha-1)x_4}dx_3dx_4+q dx_4^2\,,
\end{array}
$

\smallskip\noindent
where $a$, $b$, $c$, $q$ and $\alpha$ are arbitrary constants with $a q\neq 0$.

\medskip
\item[(A.3)] An open subset $\mathfrak{U}\subset\mathbb{R}^4$ with coordinates $(x_1,x_2,x_3,x_4)$ and metric tensor

\smallskip\noindent
$\displaystyle
g_+=2 a e^{2 x_3}\,dx_1 dx_4+a e^{2 x_3}\cos(x_4)^2dx_2^2+b dx_3^2 
+2c dx_3dx_4 +q\,dx_4^2\,, 
$

\smallskip\noindent
where  $a$, $b$, $c$ and $q$ are arbitrary constants with $ab\neq 0$ and the open set $\mathfrak{U}=\{ (x_1,x_2,x_3,x_4)\in\mathbb{R}^4;\cos(x_4)\neq 0\}$, or 

\smallskip\noindent
$\displaystyle
g_- = 2 a e^{2 x_3}\,dx_1 dx_4+a e^{2 x_3}\cosh(x_4)^2dx_2^2+b dx_3^2+2c dx_3dx_4 + q\, dx_4^2\,,
$

\smallskip\noindent
where  $a$, $b$, $c$ and $q$ are arbitrary constants with $ab\neq 0$ and $\mathfrak{U}=\mathbb{R}^4$.

\medskip
\item[(A.4)]$\mathbb{R}^4$ with coordinates $(x_1,x_2,x_3,x_4)$ and metric tensor

\smallskip\noindent
$\displaystyle
\begin{array}{l}
g=\left(\frac{a}{2}x_4^2+4bx_2^2+a\right)dx_1^2+4bx_2dx_1dx_2+ax_2(4+x_4^2)dx_1dx_3\\
\noalign{\medskip}
\phantom{g} +a(1+2x_2x_3)x_4dx_1dx_4+bdx_2^2+\frac{a}{2}(4+x_4^2)dx_2dx_3\\
\noalign{\medskip}
\phantom{g}+ax_3x_4dx_2dx_4+\frac{a}{2}dx_4^2\,,
\end{array}
$

\smallskip\noindent
where $a$ and $b$ are arbitrary constants with $a \neq 0$.

\medskip
\item[(A.5)] $(\mathbb{R}^2\setminus\{(0,0)\})\times \mathbb{R}^2$ with coordinates $(x_1,x_2,x_3,x_4)$ and metric tensor

\smallskip\noindent
$\displaystyle
\begin{array}{l}
g=-\frac{ax_4}{4x_2}dx_1dx_2+\frac{a}{4}dx_1dx_4+\frac{a(2+2x_1x_4+x_3^2)}{8x_2^2}dx_2^2\\
\noalign{\medskip}
\phantom{g} -\frac{ax_3}{4x_2}dx_2dx_3-\frac{ax_1}{4x_2}dx_2dx_4+\frac{a}{8}dx_3^2\, ,
 \end{array}
$

\smallskip\noindent
where $a\neq 0$ is an arbitrary constant.

\medskip
\item[(B.1)]$\mathbb{R}^4$ with coordinates $(x_1,x_2,x_3,x_4)$ and metric tensor

\smallskip\noindent
$\displaystyle
\begin{array}{l}
g=\left(q(x_3^2+4x_2x_3x_4+4x_2^2x_4^2)
\right.
\\
\noalign{\medskip}
\phantom{.......}
\left.
+4cx_2x_3+8cx_2^2x_4+2ax_3+4bx_2^2\right)dx_1^2\\
\noalign{\medskip}
\phantom{g}+2(q(x_3x_4+2x_2x_4^2)+4cx_2x_4+cx_3+2bx_2)dx_1dx_2\\
\noalign{\medskip}
\phantom{g}+2(q(x_3+2x_2x_4)+2cx_2+a)dx_1dx_3+4ax_2dx_1dx_4\\
\noalign{\medskip}
\phantom{g}+(qx_4^2+2cx_4+b)dx_2^2+2(qx_4+c)dx_2dx_3+2adx_2dx_4+qdx_3^2\,,
 \end{array}
$

\smallskip\noindent
where $a$, $b$, $c$ and $q$ are arbitrary constants with $a \neq 0$. 

\medskip
\item[(B.2)]$\mathfrak{U}=\{ (x_1,x_2,x_3,x_4)\in\mathbb{R}^4;x_4\neq \pm 2\}$ with coordinates $(x_1,x_2,x_3,x_4)$ and metric tensor

\smallskip\noindent
$\displaystyle
\begin{array}{l}
g=\left(a-\frac{ax_4^2}{2}+4bx_2^2\right)dx_1^2+4bx_2dx_1dx_2-a x_2(x_4^2-4)dx_1dx_3\\
 \noalign{\medskip}
\phantom{g}-a(1+2x_2x_3)x_4dx_1dx_4+bdx_2^2-\frac{1}{2}a(x_4^2-4)dx_2dx_3\\
\noalign{\medskip}
\phantom{g}-ax_3x_4dx_2dx_4-\frac{1}{2}adx_4^2\,,
 \end{array}
$

\smallskip\noindent
where $a$ and  $b$ are arbitrary constants with $a \neq 0$.

\medskip
\item[(B.3)]$\mathbb{R}^4$ with coordinates $(x_1,x_2,x_3,x_4)$ and metric tensor

\smallskip\noindent
$\displaystyle
\begin{array}{l}
g=-2ae^{-x_2}x_3dx_1dx_2+2ae^{-x_2}dx_1dx_3+2(2bx_3^2-ax_4)dx_2^2\\
\noalign{\medskip}
\phantom{g}  -4bx_3dx_2dx_3+2adx_2dx_4+bdx_3^2\,,
   \end{array}
$

\smallskip\noindent
where $a$ and $b$ are arbitrary constants with $a \neq 0$. 
\end{enumerate}
\end{theorem}

It is worth to emphasize that the spaces (A.1)-(A.3) admit metrics both of Lorentzian and neutral signature depending on the values of the constants defining the corresponding metrics. Metrics (A.4) and (A.5) are always Lorentzian, while metrics (B.1)-(B.3) are of neutral signature $(2,2)$.

Fels and Renner \cite{Fels-Renner06} classified the Einstein  non-reductive  four-dimensional  homogeneous spaces, showing that they must be of Type (A.2) or (B.3) (see also \cite{Calvaruso-Fino, Calvaruso-Zaeim14-2}). 
\begin{theorem}\label{th: Einstein}
Let $(M,g)$ be a homogeneous space of dimension four as in Theorem \ref{th:1-1}. Then $(M,g)$ is Einstein if and only if it has constant sectional curvature or it corresponds to one of the following:
\begin{enumerate}
\item Type $(A.2)$ with $\alpha=\frac{2}{3}$.
\item Type $(B.1)$ with $q=c=0\neq b$ or $q\neq 0$ and $b=\frac{c^2}{q}$.
\item Type $(B.3)$ with $b\neq 0$.
\end{enumerate}
In all the cases, the manifold is of neutral signature. 
\end{theorem}

We refer to Remark \ref{th: constant K} for a description of all non-reductive homogeneous spaces of constant sectional curvature. 
Some generalizations of the Einstein condition were studied in \cite{Calvaruso-Fino} and \cite{Calvaruso-Zaeim14} showing which of these manifolds admit Ricci solitons. 

The main goal of this paper is to study the conformal geometry of these spaces aimed to describe all the conformally Einstein non-reductive homogeneous spaces. Clearly the Einstein cases mentioned above as well as the locally conformally flat cases already described in \cite{Calvaruso-Zaeim14-2} should be discarded, since they all are conformally Einstein. 

\begin{theorem}\label{th: conformally flat}
Let $(M,g)$ be a non-reductive homogeneous space of dimension four. If $(M,g)$ is locally conformally flat then it is  of constant sectional curvature or it corresponds to one of the following cases in Theorem \ref{th:1-1}:
\begin{enumerate}
\item Type $(A.1)$ with $b=0$.
\item Type $(A.2)$ with $\alpha=2$ and $b\neq 0$.
\end{enumerate}
\end{theorem}

\subsection{The Bach Tensor}\label{se:2-2}
Let $(M^n,g)$ be a pseudo-Riemannian manifold of dimension $n$. The \emph{Schouten tensor} $\mathfrak{S}$ is the symmetric $(0,2)$-tensor field given by $\mathfrak{S}=\rho -\frac{\tau}{2(n-1)}\, g$. The failure of the Schouten tensor to be Codazzi (i.e., its covariant derivative is totally symmetric) is measured by the $(0,3)$-Cotton tensor given by  $\displaystyle\mathfrak{C}_{ijk}=(\nabla_i \mathfrak{S})_{jk} - (\nabla_j \mathfrak{S})_{ik}$.

Let $W$ denote the Weyl conformal curvature tensor and define the Ricci-contraction $W[\rho]$ to be the symmetric tensor field of type $(0,2)$ given by
$W[\rho](X,Y)$ $=$ $\sum_{ij}\varepsilon_i\varepsilon_j W(E_i,X,Y,E_j)\rho(E_i,E_j)$, where $\{ E_i\}$ is an orthonormal basis and $\varepsilon_i=g(E_i,E_i)$. 
Then the \emph{Bach tensor} is defined by
\begin{equation}
\label{eq:Bach-tensor}
\bach=\operatorname{div}_1\operatorname{div}_4 W+\frac{n-3}{n-2} W[\rho]\,,
\end{equation}
where $\operatorname{div}$ is the divergence operator. One has the coordinate description
$\displaystyle
\bach_{ij}$  $=$ $\nabla^k\nabla^\ell W_{kij\ell}$ $+$ $\frac{1}{2}\rho^{k\ell}W_{kij\ell}\,$ or, equivalently, the Bach tensor \eqref{eq:Bach-tensor} is given by 
\[
\bach_{ij}
=\frac{1}{n-2}
\left\{
\sum_{k,\alpha=1}^n g^{k\alpha} (\nabla_\alpha \mathfrak{C})_{kij} + \sum_{k,\ell=1}^n
\left( \rho_{k\ell} \sum_{\alpha,\beta=1}^n g^{k\alpha} g^{\ell\beta} W_{i\alpha j \beta}
\right) 
\right\}\,.
\]
In dimension four, the Bach tensor is  symmetric, trace-free, divergence-free and conformally invariant (i.e., if $\bar{g}=\varphi^{-2}g$, then $\bach_{\bar{g}}=\varphi^2\bach_g$). Whenever $M$ is compact, the Bach tensor is the gradient of the functional
\[
\mathcal{W}:g\mapsto \mathcal{W}(g)=\int_M \|W_g\|^2 dVg\,,
\]
and a metric is \emph{Bach flat} if it satisfies $\bach=0$. 

Bach flat metrics are critical for the functional $\mathcal{W}$, and one has that locally conformally Einstein metrics are Bach flat.
Hence, aimed to describe all the non-reductive four-dimensional homogeneous conformally Einstein metrics, one has the following 

\begin{theorem}\label{th: Bach flat}
Let $(M,g)$ be a non-reductive homogeneous space of dimension four. If $(M,g)$ is Bach flat then it is  of locally conformally flat, Einstein or it corresponds to one of the following cases in Theorem \ref{th:1-1}:
\begin{enumerate}
\item Type $(A.1)$ with $q=0$ and $b\neq 0$, or $q=-\frac{3a}{4}$ and $b\neq 0$.
\item Type $(A.2)$ with $\alpha=1$ and $b\neq 0$.
\item Type $(A.3)$ with $\epsilon=\pm 1$ and $b\neq\mp q$.
\item Type $(B.1)$ with $q=0\neq c$.
\end{enumerate}
\end{theorem}

\begin{remark}\rm
The spaces in Theorem \ref{th: Bach flat} admit Lorentzian metrics only in cases (A.2) and (A.3). All the cases in  Theorem \ref{th: Bach flat} admit metrics of neutral signature. 
\end{remark}

\begin{remark}\label{th: half conformally flat}\rm
A special class of Bach flat spaces is that of half conformally flat manifolds. While half conformally flat Lorentzian metrics are locally conformally flat, there are many strictly half conformally flat examples in the Riemannian and neutral signature settings.
Recall that a four-dimensional manifold is half conformally flat if and only it  is \emph{conformally Osserman} \cite{BGi}, i.e., the spectrum of the conformal Jacobi operators $J_W(x)(\,\cdot\,)=W(\,\cdot\,,x)x$ is constant on the unit pseudo-spheres $S^\pm(T_pM)$ at each point $p\in M$ (see \cite{Ni} and the references therein).  

An explicit calculation of the conformal Jacobi operators shows that a homogeneous four-manifold 
$(M,g)$ as in Theorem \ref{th:1-1} is half conformally flat and not locally conformally flat if and only if it corresponds to one of the following cases:
\begin{enumerate}
	\item[] Type (A.1) with $q=0\neq b$ or $q=-\frac{3}{4}a$ and $b\neq 0$. 
	\item[] Type (B.1) with $q=c=0\neq b$, or $q=0\neq c$, or $q\neq 0$ and $b=\frac{c^2}{q}$.
	\item[] Type (B.3) with $b\neq 0$.
\end{enumerate}
This agrees with the description of (anti-) self-dual non-reductive homogeneous spaces in \cite[Theorem 4.1]{Calvaruso-Zaeim14-2}. 
Moreover,
the conformal Jacobi operators are two-step nilpotent in all cases but the one corresponding to Type (B.1) with $q\neq 0$ and $b=\frac{c^2}{q}$ where they diagonalize.

It is worth to mention that in some of the cases above the manifold is also Einstein and thus \emph{pointwise Osserman}, i.e., the spectrum of the Jacobi operators $J(x)(\,\cdot\,)=R(\,\cdot\,,x)x$ is constant on the unit pseudo-spheres $S^\pm(T_pM)$ at each point $p\in M$  (see \cite{GRKVL} for further information about Osserman manifolds).

More precisely, a four-dimensional homogeneous  manifold as in Theorem \ref{th:1-1} is Osserman if and only if it is of constant sectional curvature (cf. Remark \ref{th: constant K}) or otherwise:
\begin{enumerate}
\item[(i)] $(M,g)$ is of Type (B.1) with $q=c=0\neq b$, in which case the Jacobi operators are two-step nilpotent, or
\item[(ii)] $(M,g)$ is of Type (B.1) with $q\neq 0$ and $b=\frac{c^2}{q}$. In this case, for any unit spacelike vector the corresponding Jacobi operator $J(x)(\,\cdot\,)=R(\,\cdot\,,x)x$ is diagonalizable with eigenvalues $\{0,\varepsilon_x\frac{q}{a^2},\frac{1}{4}\varepsilon_x\frac{q}{a^2}, \frac{1}{4}\varepsilon_x\frac{q}{a^2}\}$, thus locally isometric to a complex or paracomplex space form \cite{GRKVL}. A long but straightforward calculation shows that for any non-null vector $x$, the vector space $\operatorname{span}\langle x\rangle\oplus\operatorname{ker} ( J(x)-\varepsilon_x\frac{q}{a^2}\Id )$ is of Lorentzian signature. Hence $(M,g)$ is a paracomplex space form. Or
\item[(iii)] $(M,g)$ is of Type (B.3) with $b\neq 0$, in which case the Jacobi operators are two-step nilpotent. 
\end{enumerate}
Moreover, it is worth to emphasize that in all the cases above the manifold is locally symmetric.
\end{remark}

\subsection{Conformally Einstein manifolds}\label{se:2-3}

A central problem in conformal geometry is to decide whether a given manifold $(M,g)$ is in the conformal class of an Einstein manifold, i.e., if there exists a (locally defined) smooth function $\varphi$ such that $\bar{g}=\varphi^{-2}g$ is Einstein. It is well-known that any conformal transformation  preserves the Weyl tensor of type $(1,3)$, but neither the connection nor the curvature tensor remain invariant. The Ricci tensor ${}^{g}\rho$ changes under a conformal transformation $\bar{g}=\varphi^{-2}g$ as 
\begin{equation}
{}^{\bar{g}}\rho-{}^{g}\rho=\varphi^{-2}\left((n-2)\varphi \Hes_\varphi+\left( \varphi\Delta\varphi-(n-1)\|\nabla\varphi\|^2\right)g \right)\,,
\end{equation}
where $\Hes_\varphi=\nabla d\varphi$ is the Hessian of $\varphi$, $\Delta=\operatorname{trace}\Hes$ is the Laplacian and $\nabla\varphi$ is the gradient of $\varphi$.

A conformal transformation $\varphi$ maps an Einstein metric $g$ to another Einstein metric $\bar{g}=\varphi^{-2}g$ if and only if
\begin{equation}
\Hes_{\varphi}=\frac{\Delta\varphi}{n}g\,.
\end{equation} 
This equation was studied by Brinkmann showing that if $g(\nabla\varphi,\nabla\varphi)\neq0$ then the manifold is locally a warped product and if $\nabla\varphi$ is a null vector field then it is parallel and $(M,g)$ is a Walker manifold \cite{Brinkmann25}, \cite{Walker}.

$(M,g)$ is said to be \emph{(locally) conformally Einstein} if every point $p\in M$
has an open neighborhood $\mathcal{U}$ and a positive smooth function $\varphi$ defined on $\mathcal{U}$ such that $(\mathcal{U},\bar{g} =
\varphi^{-2} g)$ is Einstein.
Brinkmann \cite{Brinkmann24} showed that a manifold is conformally Einstein if and only if
the equation
\begin{equation}\label{eq: conformal Einstein equation}
(n-2)\Hes_\varphi +\varphi\,\rho= \frac{1}{n}\{(n-2)\Delta\varphi+\varphi\,\tau\}g
\end{equation}
has a non-constant solution.

It is important to emphasize that although any locally conformally Einstein metric is Bach flat, there are examples of strictly Bach flat manifolds, i.e., they are neither half conformally flat nor locally conformally Einstein (see, for example \cite{AGS, nuevo,LNu} and references therein).
Now, our main result may be stated as follows.

\bigskip
\noindent\textbf{Main Theorem.}\label{th:main}\emph{
Let $(M,g)$ be a four dimensional non-reductive homogeneous space. $(M,g)$ is in the conformal class of an Einstein manifold if and only if $(M,g)$ is Einstein, locally conformally flat, or is locally isometric to one of this spaces:
\begin{enumerate}
\item Type $(A.1)$ with $q=0$ and $b\neq 0$, or $q=-\frac{3a}{4}$ and $b\neq 0$.
\item Type $(A.2)$ with $\alpha=1$ and $b\neq 0$.
\item Type $(A.3)$ with $\epsilon=\pm 1$ and $b\neq\mp q$.
\end{enumerate}
Moreover, all the cases $(1)-(3)$ are in the conformal class of a Ricci flat metric which is unique (up to a constant) only in Type (A.1) with $q=0$. Otherwise the space of Ricci flat conformal metrics is either two or three-dimensional.}

\section{Curvature of non-reductive four-dimensional homogeneous spaces}\label{se:3}

In this section we briefly schedule some basic facts about the curvature of non-reductive homogeneous spaces. All the curvature expressions are obtained after some straightforward calculations that we omit. We consider separately all the possibilities in Theorem \ref{th:1-1} and analyze the Ricci, the Cotton, the Weyl and the Bach tensor case by case. As a consequence, one obtains the proofs of Theorem \ref{th: Einstein}, Theorem \ref{th: conformally flat} and Theorem \ref{th: Bach flat}.

\subsection{Non-reductive homogeneous manifolds admitting Lorentzian and neutral signature metrics.} With the notation of  Theorem~\ref{th:1-1} at hand, the non-reductive four-dimensional homogeneous manifolds admitting both Lorentzian and neutral signature metrics are those corresponding to types (A.1), (A.2) and (A.3).

\subsubsection{\textbf{\rm Type (A.1)}} Consider the metric tensor
\begin{equation}
\label{eq:A.1-1}
\begin{array}{l}
g=(4b x_2^2+a)\,dx_1^2+4bx_2\,dx_1dx_2-(4ax_2x_4-4cx_2+a)\,dx_1dx_3 \\
\noalign{\medskip}
\phantom{g=}+4ax_2\,dx_1dx_4+b\,dx_2^2-2(ax_4-c)\,dx_2dx_3+2a\,dx_2dx_4+q\,dx_3^2\,.
\end{array}
\end{equation}
It immediately follows from the above expression that $\operatorname{det}(g)=\frac{1}{4}a^3(a-4q)$, which shows that the metric \eqref{eq:A.1-1} is Lorentzian if $a(a-4q)<0$ and of neutral signature otherwise. Further observe that the restriction $a(a-4q)\neq 0$ in Theorem \ref{th:1-1}-(A.1) ensures that $g$ is non-degenerate.

The 
Ricci operator is given by
\begin{equation}\label{eq: Ricciop A.1}
\Ric =\frac{1}{a}
\left(\begin{array}{cccc}
-2 &0&1&0\\
\noalign{\medskip}
0&-2&-2 x_2&0\\
\noalign{\medskip}
0&0&0&0\\
\noalign{\medskip}
\frac{8b(a+4q)x_2}{a(a-4q)}&\frac{4b(a+4q)}{a(a-4q)}&\frac{2(ax_4-c)}{a}&-2
\end{array}\right)\,,
\end{equation}
from where it follows that $(M,g)$ is not Einstein.
The non-zero components of the Cotton tensor are given by (up to symmetries):
\begin{equation}\label{eq: Cotton A.1}
\begin{array}{l}
\mathfrak{C}_{121}=-\frac{24 b x_2 (a+4 q)}{a (a-4 q)}\,,\quad \mathfrak{C}_{122}=-\frac{12 b (a+4 q)}{a (a-4 q)}\,,\quad \mathfrak{C}_{131}=\frac{64 b q x_2^2}{a^2-4 a q}\,,\\ 
\noalign{\bigskip}
\mathfrak{C}_{132}=\mathfrak{C}_{231}=\frac{32 b q x_2}{a^2-4 a q}\,,\quad
\mathfrak{C}_{232}=\frac{16 b q}{a^2-4 a q}\,,
\end{array}
\end{equation}
and the non-zero components of the Weyl tensor are (up the usual symmetries):
\begin{equation}\label{eq: Weyl A.1}
\begin{array}{lll}
W_{1212}=\frac{8 b q-6 a b}{a-4 q}\,,&\! W_{1213}=-\frac{16 b q x_2}{a-4 q}\,,& W_{1223}=-\frac{8 b q}{a-4 q}\,,\\
\noalign{\medskip}
W_{1313}=-\frac{8 b q x_2^2 (a+4 q)}{a (a-4 q)}\,,&\! W_{1323}=-\frac{4 b q x_2 (a+4 q)}{a (a-4 q)}\,,& W_{2323}=-\frac{2 b q (a+4 q)}{a (a-4 q)}\,.
\end{array}
\end{equation}
Hence, the Bach tensor is given by
\begin{equation}\label{eq: Bach A.1}
\mathfrak{B}=\left(
\begin{array}{cccc}
 -\frac{256 b\, q \,(3 a+4 q)\,x_2^2 }{a^2 (a-4 q)^2} & -\frac{128 b\, q\, (3 a+4 q) x_2}{a^2 (a-4 q)^2} & 0 & 0 \\
 -\frac{128 b\, q\, (3 a+4 q)\, x_2 }{a^2 (a-4 q)^2} & -\frac{64 b\, q\, (3 a+4 q)}{a^2 (a-4 q)^2} & 0 & 0 \\
 0 & 0 & 0 & 0 \\
 0 & 0 & 0 & 0
\end{array}
\right)\,.
\end{equation}

An immediate consequence of previous expression is that \emph{a Type $(A.1)$ homogeneous space is Bach flat if and only if one of the following constraints holds: $b=0$, $q=0$ or $q=-\frac{3a}{4}$. Moreover:
\begin{enumerate}
\item If $b=0$, then \eqref{eq: Weyl A.1} shows that $(M,g)$ is locally conformally flat.
\item If $b\neq 0$, then $(M,g)$ is neither locally conformally flat nor Einstein.
\end{enumerate}
}

\subsubsection{\textbf{\rm Type (A.2)}} Consider the metric tensor
\begin{equation}\label{eq:A.2-1}
\begin{array}{l}
g =  -2 a e^{2\alpha x_4}\,dx_1 dx_3+a e^{2\alpha x_4}dx_2^2+b e^{2(\alpha-1)x_4}dx_3^2 
\\
\noalign{\medskip}
\phantom{g=}
+2c e^{(\alpha-1)x_4}dx_3dx_4+q dx_4^2\,.
\end{array}
\end{equation}
It immediately follows from the above expression that $\operatorname{det}(g)=-a^3 q\, e^{6\alpha x_4}$, which shows that the metric \eqref{eq:A.2-1} is Lorentzian if $aq>0$ and of neutral signature otherwise. Further observe that the restriction $aq\neq 0$ in Theorem \ref{th:1-1}-(A.2) ensures that $g$ is non-degenerate.

The  
Ricci operator is given by
\begin{equation}\label{eq: Ricciop A.2}
\Ric=-\frac{3\alpha^2}{q}\left(\begin{array}{cccc}
1&0&\frac{b(3\alpha-2)}{3a\alpha^2}e^{-2x_4}&0\\
\noalign{\medskip}
0&1&0&0\\
\noalign{\medskip}
0&0&1&0\\
\noalign{\medskip}
0&0&0&1
\end{array}\right)\,.
\end{equation}
Hence $(M,g)$ is Einstein if and only if $b=0$ (with scalar curvature $\tau=-12\frac{\alpha^2}{q}$), or $\alpha=\frac{2}{3}$ (with scalar curvature $\tau=-\frac{16}{3q}$), or Ricci flat if $\alpha=0$.

The only non-zero component of the Cotton tensor is given by
\begin{equation}\label{eq: Cotton A.2}
\mathfrak{C}_{343}=-\frac{(\alpha-2)(3\alpha-2)\, b e^{2 (\alpha -1) x_4}}{q}\,,
\end{equation}
and the non-zero components of the Weyl tensor are given by
\begin{equation}\label{eq: Weyl A.2}
W_{2323}=-\frac{(\alpha -2) a\, b e^{2 (2 \alpha -1) x_4}}{2 q}\,,\quad W_{3434}=\frac{1}{2} (\alpha -2) b e^{2 (\alpha -1) x_4}\,.
\end{equation}
Finally the Bach tensor is expressed with respect to the coordinate basis as
\begin{equation}\label{eq:bach A.2}
\mathfrak{B}=\left(
\begin{array}{cccc}
 0 & 0 & 0 & 0 \\
 0 & 0 & 0 & 0 \\
 0 & 0 & \frac{(\alpha-2) (\alpha-1) (3 \alpha-2) b e^{2 (\alpha-1) x_4}}{q^2} & 0 \\
 0 & 0 & 0 & 0
\end{array}
\right)\,.
\end{equation}

Hence \emph{a homogeneous space of Type $(A.2)$ is Bach flat if and only if $b=0$, $\alpha=\frac{2}{3}$, $\alpha=1$ or $\alpha=2$. Moreover:
\begin{enumerate}
\item If $b=0$ then \eqref{eq: Ricciop A.2} and \eqref{eq: Weyl A.2} show that the manifold is of constant sectional curvature $K=-\frac{\alpha^2}{q}$. 
\item If  $\alpha=\frac{2}{3}$, then $W_{3434}=-\frac{2}{3} b e^{-\frac{2 x_4}{3}}$ and hence the manifold is not locally conformally flat, unless  $b=0$. 
\item If $\alpha=1$ then $W_{3434}=-\frac{b}{2}$, which shows that $(M,g)$ is not locally conformally flat unless $b= 0$. 
\item If $\alpha=2$, then \eqref{eq: Weyl A.2} shows that $(M,g)$ is locally conformally flat but not Einstein unless $b=0$.
\end{enumerate}
}

\subsubsection{\textbf{\rm Type (A.3)}} Two distinct cases have to be considered for Type (A.3) metrics.
Let $\mathfrak{U}$ be the open set in $\mathbb{R}^4$ determined by 
$\mathfrak{U}=\{ (x_1,x_2,x_3,x_4)\in\mathbb{R}^4;\cos(x_4)\neq 0\}$ and the metric tensor
\begin{equation}\label{eq:A.3-1}
g_+=2 a e^{2 x_3}\,dx_1 dx_4+a e^{2 x_3}\cos(x_4)^2dx_2^2+b dx_3^2 
+2c dx_3dx_4 +q\,dx_4^2\,.
\end{equation}
Now $\operatorname{det}(g_+)=-a^3 b\cos(x_4)^2\, e^{6x_3}$ shows that the metric \eqref{eq:A.3-1} is Lorentzian if $ab>0$ and of neutral signature otherwise. Further observe that the restriction $ab\neq 0$ in Theorem \ref{th:1-1}-(A.3) ensures that $g_+$ is non-degenerate.

The 
Ricci operator is given by
\begin{equation}\label{eq: Ricciop A.31}
\Ric=-\frac{3}{b}\left(
\begin{array}{cccc}
 1 & 0 & 0 & -\frac{(b+q)e^{-2x_3}}{3a} \\
 0 & 1 & 0 & 0 \\
 0 & 0 & 1 &0 \\
 0 & 0 &0 & 1
\end{array}
\right)\,,
\end{equation}
and thus $(M,g)$ is Einstein if and only if $b=-q$.

The only non-zero component of the Cotton tensor is
\begin{equation}\label{eq: Cotton A.31}
\mathfrak{C}_{344}=-\frac{b+q}{b}\,,
\end{equation}
and the non-zero components of the Weyl tensor are:
\begin{equation}\label{eq: Weyl A.31}
W_{2424}=\frac{a e^{2 x_3} (b+q) \cos\left(x_4\right)^2}{2 b}\,,\quad W_{3434}=-\frac{b + q}{2}\,.
\end{equation}

Now, a long but straightforward computation shows that 
\emph{$(M,g_+)$ is always Bach flat. Moreover, $(M,g_+)$ is locally conformally flat if and only if $b=-q$ by Equation \eqref{eq: Weyl A.31}, in which case it is Einstein and thus of constant sectional curvature $K=\frac{1}{q}$.}

\bigskip

Now we consider the second case for Type (A.3) metrics.
Let $M$ be $\mathbb{R}^4$ with metric tensor
\begin{equation}\label{eq:A.3-2}
g_- = 2 a e^{2 x_3}\,dx_1 dx_4+a e^{2 x_3}\cosh(x_4)^2dx_2^2+b dx_3^2+2c dx_3dx_4 + q\, dx_4^2 \,.
\end{equation}
Next  $\operatorname{det}(g_-)=-a^3 b\cosh(x_4)^2\, e^{6x_3}$ shows that the metric \eqref{eq:A.3-2} is Lorentzian if $ab>0$ and of neutral signature otherwise. Further observe that the restriction $ab\neq 0$ in Theorem \ref{th:1-1}-(A.3) ensures that $g_-$ is non-degenerate.

The 
Ricci operator is given by
\begin{equation}\label{eq: Ricciop A.3-1}
\Ric=-\frac{3}{b}\left(
\begin{array}{cccc}
 1 & 0 & 0 & \frac{(b-q)e^{-2x_3}}{3a} \\
 0 & 1 & 0 & 0 \\
 0 & 0 & 1 &0 \\
 0 & 0 &0 & 1
\end{array}
\right)\,,
\end{equation}
and thus $(M,g)$ is Einstein if and only if $b=q$.

The only non-zero component of the Cotton tensor is
\begin{equation}\label{eq: Cotton A.3-1}
\mathfrak{C}_{344}=1-\frac{q}{b}\,,
\end{equation}
and the non-vanishing components of the Weyl tensor are:
\begin{equation}\label{eq: Weyl A.3 -1}
W_{2424}=-\frac{a e^{2 x_3} (b-q) \cosh\left(x_4\right)^2}{2 b}\,,\quad W_{3434}=\frac{b-q}{2}\,.
\end{equation}
Furthermore, a long but straightforward computation shows that 
\emph{$(M,g_-)$ is always Bach flat. Moreover, $(M,g_-)$ is locally conformally flat if and only if $b=q$ by Equation \eqref{eq: Weyl A.3 -1}, in which case it is Einstein and thus of constant sectional curvature $K=-\frac{1}{q}$.}

\emph{Hence any Einstein Type (A.3) manifold is necessarily of constant sectional curvature.}

\subsection{Non-reductive homogeneous spaces admitting only Lorentzian metrics}

\subsubsection{\textbf{\rm Type (A.4)}} Consider the metric tensor
\begin{equation}\label{eq:A.4-1}
\begin{array}{l}
g=\left(\frac{a}{2}x_4^2+4bx_2^2+a\right)dx_1^2+4bx_2dx_1dx_2+ax_2(4+x_4^2)dx_1dx_3\\
\noalign{\medskip}
\phantom{g=} +a(1+2x_2x_3)x_4dx_1dx_4+bdx_2^2+\frac{a}{2}(4+x_4^2)dx_2dx_3\\
\noalign{\medskip}
\phantom{g=}+ax_3x_4dx_2dx_4+\frac{a}{2}dx_4^2\,.
\end{array}
\end{equation}
It   follows from the above expression that $\operatorname{det}(g)=-\frac{1}{32}a^4 (4+x_4^2)^2$, which shows that the metric \eqref{eq:A.4-1} is Lorentzian. Further observe that the restriction $a\neq 0$ in Theorem \ref{th:1-1}-(A.4) ensures that $g$ is non-degenerate.
 
The 
Ricci operator is given by
\begin{equation}\label{eq: Ricciop A.4}
\Ric=-\frac{3}{a}  \left(\begin{array}{cccc}
1&0&0&0\\
\noalign{\medskip}
0&1&0&0\\
\noalign{\medskip}
\frac{40 b x_2}{3a(x_4^2+4)} & \frac{20 b}{3a(x_4^2+4)} &1&0\\
\noalign{\medskip}
0&0&0&1
\end{array}\right)\,,
\end{equation}
which shows that $(M,g)$ is Einstein if and only if $b=0$.

The only non-zero components of the Cotton tensor are given by
\begin{equation}
 \mathfrak{C}_{121}=\frac{30 b x_2}{a}\,,\qquad \mathfrak{C}_{122}=\frac{15 b}{a}\,,
\end{equation}
and the non-zero components of the Weyl tensor are
\begin{equation}\label{eq: Weyl A.4}
\begin{array}{lll}
W_{1212}=\frac{3}{4} b \left(x_4^2-2\right)\,,& W_{1214}=-\frac{3}{2} b x_2 x_4\,, & W_{1224}=-\frac{3 b x_4}{4}\,,\\
\noalign{\medskip}
W_{1414}=3 b x_2^2\,,& W_{1424}=\frac{3 b x_2}{2}\,,& W_{2424}=\frac{3 b}{4}\,.
\end{array}
\end{equation}
The Bach tensor is given by
\begin{equation}\label{eq: Bach A.4}
\mathfrak{B}=\left(
\begin{array}{cccc}
 -\frac{120 b x_2^2}{a^2} & -\frac{60 b x_2}{a^2} & 0 & 0 \\
 -\frac{60 b x_2}{a^2} & -\frac{30 b}{a^2} & 0 & 0 \\
 0 & 0 & 0 & 0 \\
 0 & 0 & 0 & 0
\end{array}
\right)\,.
\end{equation}
Hence, \emph{a Type $(A.4)$ metric is Bach flat if and only if $b=0$, in which case $(M,g)$ is locally conformally flat by \eqref{eq: Weyl A.4}, and thus of constant sectional curvature $K=-\frac{1}{a}$}, as the Ricci operator shows.

\subsubsection{\textbf{\rm Type (A.5)}} Let $M=(\mathbb{R}^2\setminus\{(0,0)\})\times \mathbb{R}^2$ and let $(x_1,x_2,x_3,x_4)$ be the coordinates. Consider the  metric tensor
\begin{equation}\label{eq:A.5-1}
\begin{array}{l}
g=-\frac{ax_4}{4x_2}dx_1dx_2+\frac{a}{4}dx_1dx_4+\frac{a(2+2x_1x_4+x_3^2)}{8x_2^2}dx_2^2\\
\noalign{\medskip}
\phantom{g=} -\frac{ax_3}{4x_2}dx_2dx_3-\frac{ax_1}{4x_2}dx_2dx_4+\frac{a}{8}dx_3^2\,.
 \end{array}
\end{equation}
Since $\operatorname{det}(g)=-\frac{a^4}{2048x_2^2}$, the metric \eqref{eq:A.5-1} is Lorentzian and the restriction $a\neq 0$ in Theorem \ref{th:1-1}-(A.5) ensures that $g$ is non-degenerate.
 
The Ricci tensor is given by
\begin{equation}\label{eq: Ricci A.5}
\rho=\left(
\begin{array}{cccc}
 0 & \frac{3 x_4}{2 x_2} & 0 & -\frac{3}{2} \\
 \frac{3 x_4}{2 x_2} & -\frac{3 \left(x_3^2+2 x_1 x_4+2\right)}{2 x_2^2} & \frac{3 x_3}{2 x_2} & \frac{3 x_1}{2 x_2} \\
 0 & \frac{3 x_3}{2 x_2} & -\frac{3}{2} & 0 \\
 -\frac{3}{2} & \frac{3 x_1}{2 x_2} & 0 & 0
\end{array}
\right)\,,
\end{equation} 
from where it follows that the corresponding Ricci operator is a multiple of the identity,
$\Ric=-\frac{12}{a}\operatorname{Id}$, and thus Einstein. Moreover, the Weyl tensor vanishes identically and therefore, any \emph{Type (A.5) metric is always of constant sectional curvature $K=-\frac{4}{a}$}.

\subsection{Non-reductive  four-dimensional homogeneous manifolds admitting only neutral signature metrics.}
There exist three different families of non-reductive homogeneous four-manifolds which admit exclusively neutral signature metrics. 

\subsubsection{\textbf{\rm Type (B.1)}} Let $M=\mathbb{R}^4$ with coordinates $(x_1,x_2,x_3,x_4)$ and metric tensor
\begin{equation}\label{eq:B.1-1}
\begin{array}{l}
g\!=\!\left(q(x_3^2+4x_2x_3x_4+4x_2^2x_4^2)+4cx_2x_3+8cx_2^2x_4+2ax_3+4bx_2^2\right)dx_1^2\!\!\\
\noalign{\medskip}
\phantom{g=}+2(q(x_3x_4+2x_2x_4^2)+4cx_2x_4+cx_3+2bx_2)dx_1dx_2\\
\noalign{\medskip}
\phantom{g=}+2(q(x_3+2x_2x_4)+2cx_2+a)dx_1dx_3+4ax_2dx_1dx_4\\
\noalign{\medskip}
\phantom{g=}+(qx_4^2+2cx_4+b)dx_2^2+2(qx_4+c)dx_2dx_3+2adx_2dx_4+qdx_3^2\,.
 \end{array}
\end{equation}
Since $\operatorname{det}(g)=a^4$ and the component $g_{44}=0$, the metric \eqref{eq:B.1-1} is of neutral signature and the restriction $a\neq 0$ in Theorem \ref{th:1-1}-(B.1) ensures that $g$ is non-degenerate.
 
The Ricci operator is given by
\begin{equation}\label{eq: Ricciop B.1}
\Ric=\left(\begin{array}{cccc}
\frac{3q}{2a^2} &0&0&0\\
\noalign{\medskip}
0&\frac{3q}{2a^2} &0&0\\
\noalign{\medskip}
0&0&\frac{3q}{2a^2}&0\\
\noalign{\medskip}
\frac{15}{a^3}x_2(bq-c^2) & \frac{15}{2 a^3}(bq-c^2) &0&\frac{3q}{2a^2}
\end{array}\right)\,,
\end{equation}
from where it follows that $(M,g)$ is Einstein if and only if $c^2-bq=0$.

The non-zero components of the Cotton tensor are given by
\begin{equation}\label{eq: Cotton B.1}
\begin{array}{l}
\mathfrak{C}_{121}=\frac{15 x_2 \left(6 a-q x_3\right) \left(c^2-b q\right)}{2 a^3}\,,\quad \mathfrak{C}_{122}=\frac{15 \left(6 a-q x_3\right) \left(c^2-b q\right)}{4 a^3}\,,\\
\noalign{\medskip} 
 \mathfrak{C}_{232}=\frac{15 q \left(c^2-b q\right)}{4 a^3}\,,
 \quad
 \mathfrak{C}_{131}=4 x_2^2 \,\mathfrak{C}_{232}\,,
 \quad
 \mathfrak{C}_{132}=\mathfrak{C}_{231}=2 x_2 \mathfrak{C}_{232}\,, 
\end{array}
\end{equation}
and the non-zero components of the Weyl tensor are given by
\begin{equation}\label{eq: Weyl B.1}
\begin{array}{l}
W_{1212}=\frac{-6 a^2 \left(b+2 c x_4+q x_4^2\right)+a x_3 \left(-7 b q+6 c^2-q x_4 \left(2 c+q x_4\right)\right)+5 q x_3^2 \left(c^2-b q\right)}{2 a^2}\,,\\
\noalign{\medskip}
W_{1213}=\frac{2 x_2 \left(a \left(7 b q-6 c^2+q x_4 \left(2 c+q x_4\right)\right)+10 q x_3 \left(b q-c^2\right)\right)-a \left(6 a+q x_3\right) \left(c+q x_4\right)}{4 a^2}\,,\\
\noalign{\medskip}
W_{1214}=-\frac{2 x_2 \left(6 a+q x_3\right) \left(c+q x_4\right)+q x_3 \left(2 a+q x_3\right)}{4 a}\,,\\
\noalign{\medskip}
W_{1223}=\frac{a \left(7 b q-6 c^2+q x_4 \left(2 c+q x_4\right)\right)+10 q x_3 \left(b q-c^2\right)}{4 a^2}\,,\\
\noalign{\medskip}
W_{1224}=-\frac{\left(6 a+q x_3\right) \left(c+q x_4\right)}{4 a}\,,\quad W_{1234}=-\frac{q \left(2 a+q x_3\right)}{4 a}\,,\\
\noalign{\medskip}
W_{1313}=\frac{q \left(-a^2+2 a x_2 \left(c+q x_4\right)+20 x_2^2 \left(c^2-b q\right)\right)}{2 a^2}\,,\\
\noalign{\medskip}
W_{1314}=\frac{q x_2 \left(-2 a+2 x_2 \left(c+q x_4\right)+q x_3\right)}{2 a}\,,\\
\noalign{\medskip}
W_{1323}=\frac{q \left(a \left(c+q x_4\right)+20 x_2 \left(c^2-b q\right)\right)}{4 a^2}\,,
\quad W_{1324}=\frac{q \left(x_2 \left(c+q x_4\right)-a\right)}{2 a}\,,\\
\noalign{\medskip}
W_{1334}=\frac{q^2 x_2}{2 a},\quad W_{1423}=\frac{q \left(2 x_2 \left(c+q x_4\right)+q x_3\right)}{4 a},\quad W_{1424}=-q x_2\,,\\
\noalign{\medskip}
W_{1414}=-2 q x_2^2\,,\quad W_{2334}=\frac{q^2}{4 a}\,, \quad W_{2424}=-\frac{q}{2}\,,\\
\noalign{\medskip}
W_{2323}=\frac{5 q \left(c^2-b q\right)}{2 a^2}\,,\quad W_{2324}=\frac{q \left(c+q x_4\right)}{4 a}\,.
\end{array}
\end{equation}
Hence, the Bach tensor is given by
\begin{equation}\label{eq: Bach B.1}
\mathfrak{B}=\left(
\begin{array}{cccc}
 \frac{240 q \left(c^2-b q\right) x_2^2}{a^4} & \frac{120 q \left(c^2-b q\right) x_2}{a^4} & 0 & 0 \\
 \frac{120 q \left(c^2-b q\right) x_2}{a^4} & \frac{60 q \left(c^2-b q\right)}{a^4} & 0 & 0 \\
 0 & 0 & 0 & 0 \\
 0 & 0 & 0 & 0
\end{array}
\right)\,.
\end{equation}
Thus, \emph{a Type $(B.1)$ metric is Bach flat if and only if $q=0$ or $c^2-bq=0$, in the later case being Einstein. Moreover},
\begin{enumerate}
\item \emph{If $q=0$, then the Ricci operator \eqref{eq: Ricciop B.1} is either zero or two-step nilpotent, Equation \eqref{eq: Weyl B.1} gives $W_{1224}=-\frac{3}{2}c$, thus distinguishing  the following two cases:}
\begin{enumerate}
\item \emph{If $q=0$ and $c=0$, then $(M,g)$ is Ricci flat and the only non-zero component of the Weyl tensor is $W_{1212}=-3b$. Therefore  $(M,g)$ is flat if $q=c=b=0$}. 

\emph{Otherwise, if $q=c=0\neq b$, then the Jacobi operators are two-step nilpotent. Hence $(M,g)$ is  Osserman and thus half conformally flat}.

\item \emph{If $q=0$ and $c\neq 0$, then $(M,g)$ is not locally conformally flat. Moreover the conformal Jacobi operators are nilpotent and $(M,g)$ is half conformally flat.}
\end{enumerate}
\item \emph{If $q\neq 0$ and $b=\frac{c^2}{q}$, then \eqref{eq: Weyl B.1} shows that $W_{1334}=\frac{q^2 x_2}{2 a }$ and hence $(M,g)$ is not locally conformally flat}. 
\emph{Equation \eqref{eq: Ricciop B.1} shows that $(M,g)$ is Einstein and moreover the Jacobi operator $J(x)(\,\cdot\,)=R(\,\cdot\,,x)x$ associated to any unit vector $x$ has constant eigenvalues $\{ 0,\varepsilon_x\frac{q}{a^2},\varepsilon_x\frac{q}{4a^2},\varepsilon_x\frac{q}{4a^2}\}$}. 

\emph{Moreover $(M,g)$ is locally isometric to a paracomplex space form of constant paraholomorphic sectional curvature $H=-\frac{q}{a^2}$, and thus a modified Riemannian extension as in \cite{CLGRGVL}}.
\end{enumerate}

\subsubsection{\textbf{\rm Type (B.2)}} Let 
$\mathfrak{U}=\{ (x_1,x_2,x_3,x_4)\in\mathbb{R}^4;x_4\neq \pm 2\}$ with coordinates $(x_1,x_2,x_3,x_4)$ and metric tensor
\begin{equation}\label{eq:B.2-1}
\begin{array}{l}
g=\left(a-\frac{ax_4^2}{2}+4bx_2^2\right)dx_1^2+4bx_2dx_1dx_2-a x_2(x_4^2-4)dx_1dx_3\\
 \noalign{\medskip}
\phantom{g=}-a(1+2x_2x_3)x_4dx_1dx_4+bdx_2^2-\frac{1}{2}a(x_4^2-4)dx_2dx_3\\
\noalign{\medskip}
\phantom{g=}-ax_3x_4dx_2dx_4-\frac{1}{2}adx_4^2\,.
 \end{array}
\end{equation}
Since $\operatorname{det}(g)=\frac{1}{32}a^4(x_4^2-4)^2$ and the component $g_{33}=0$, the metric \eqref{eq:B.2-1} is of neutral signature and the restriction $a\neq 0$, $x_4\neq \pm 2$  in Theorem \ref{th:1-1}-(B.2) ensures that $g$ is non-degenerate.
 
The Ricci operator is given by
\begin{equation}\label{eq: Ricciop B.2}
\Ric=-\frac{3}{a}\left(
\begin{array}{cccc}
1 &0&0&0\\
\noalign{\medskip}
0&1&0&0\\
\noalign{\medskip}
-\frac{40 b x_2}{3a(x_4^2-4)} & -\frac{20 b}{3a(x_4^2-4)} & 1&0\\
\noalign{\medskip}
0&0&0&1
\end{array}
\right)\,,
\end{equation}
which shows that $(M,g)$ is Einstein if and only if $b=0$.

The non-zero components of the Cotton tensor are given by
\begin{equation}
\mathfrak{C}_{121}=\frac{30 b x_2}{a}\,,\quad \mathfrak{C}_{122}=\frac{15 b}{a}\,,
\end{equation}
and  the non-zero components of the Weyl tensor are determined by
\begin{equation}\label{eq: Weyl B.2}
\begin{array}{lll}
W_{1212}=-\frac{3}{4} b \left(x_4^2+2\right)\,,& W_{1214}=\frac{3}{2} b x_2 x_4\,,& W_{1224}=\frac{3 b x_4}{4}\,,\\
\noalign{\medskip}
W_{1414}=-3 b x_2^2\,,& W_{1424}=-\frac{3 b x_2}{2}\,,& W_{2424}=-\frac{3 b}{4}\,.
\end{array}
\end{equation}
Hence, the Bach tensor is given by
\begin{equation}\label{eq: Bach B.2}
\mathfrak{B}=\left(
\begin{array}{cccc}
 -\frac{120 b x_2^2}{a^2} & -\frac{60 b x_2}{a^2} & 0 & 0 \\
 -\frac{60 b x_2}{a^2} & -\frac{30 b}{a^2} & 0 & 0 \\
 0 & 0 & 0 & 0 \\
 0 & 0 & 0 & 0
\end{array}
\right)\,.
\end{equation}

Now it follows from the previous expressions that \emph{a metric \eqref{eq:B.2-1} is Bach flat if and only if $b=0$, in which case it is Einstein and locally conformally flat, and thus  of constant sectional curvature $K=-\frac{1}{a}$}.

\subsubsection{\textbf{\rm Type (B.3)}} Let 
$\mathbb{R}^4$ with coordinates $(x_1,x_2,x_3,x_4)$ and metric tensor
\begin{equation}\label{eq:B.3-1}
\begin{array}{l}
g=-2ae^{-x_2}x_3dx_1dx_2+2ae^{-x_2}dx_1dx_3+2(2bx_3^2-ax_4)dx_2^2\\
\noalign{\medskip}
\phantom{g=}  -4bx_3dx_2dx_3+2adx_2dx_4+bdx_3^2\,.
   \end{array}
\end{equation}
Since $\operatorname{det}(g)=a^4 e^{-2x_2}$ and the component $g_{44}=0$, the metric \eqref{eq:B.3-1} is of neutral signature and the restriction $a\neq 0$  in Theorem \ref{th:1-1}-(B.3) ensures that $g$ is non-degenerate.
 
A straightforward calculation shows that \emph{the Ricci operator of any Type $(B.3)$ metric vanishes identically and hence they are all Ricci flat}. Thus the Cotton and the Bach tensor are also zero. The Weyl tensor is, however, not necessarily zero and the only non-zero component of the Weyl tensor is given by 
\begin{equation}
W_{2323}=-3b\,,
\end{equation}
which shows that $(M,g)$ is flat if and only if $b=0$.

\emph{Any metric of Type (B.3) with $b\neq 0$ has two-step nilpotent Jacobi operators, and thus it is Osserman}.

\begin{remark}
	\rm\label{th: constant K}
As a consequence of the expressions of the Ricci and the Weyl tensor in this section, a  four-dimensional homogeneous space $(M,g)$ as in Theorem \ref{th:1-1}  is of constant sectional curvature $K$ if and only if it corresponds to one of the following (see also \cite{Calvaruso-Fino, Calvaruso-Zaeim14-2, Fels-Renner06}):
\begin{enumerate}
	\item[] Type (A.2) with $b=0$, in which case $K=-\frac{\alpha^2}{q}$.
	\item[] Type (A.3) with $b=-\epsilon q$, in which case $K=\epsilon\frac{1}{q}$.
	\item[] Type (A.4) with $b=0$, in which case $K=-\frac{1}{a}$.
	\item[] Type (A.5), in which case $K=-\frac{4}{a}$.
	\item[] Type (B.1) with $q=c=b=0$, in which case is flat.
	\item[] Type (B.2) with $b=0$, in which case $K=-\frac{1}{a}$.
	\item[] Type (B.3) with $b=0$, in which case is flat.
\end{enumerate}
Moreover, a long but straightforward calculation shows that a four-dimensional  homogeneous space given by Theorem \ref{th:1-1}, of non-constant sectional curvature,  is locally symmetric if and only if it it is
\begin{enumerate}
\item[] Type (A.1) with $b=0$, in which case $(M,g)$ is locally conformally flat with diagonalizable Ricci operator. Hence locally isometric  to a product $\mathbb{R}\times N$, where $N$ is of constant sectional curvature $K_N=-\frac{1}{a}$, 
\end{enumerate}
or it corresponds to one of the following cases:
\begin{enumerate}
\item[] Type (B.1) with $q=c=0\neq b$, in which case $(M,g)$ is Osserman with two-step nilpotent Jacobi operators.
\item[] Type (B.1) with $q\neq 0$ and $b=\frac{c^2}{q}$, in which case $(M,g)$ is a paracomplex space form.
\item[] Type (B.3) with $b\neq 0$, in which case $(M,g)$ is Osserman with two-step nilpotent Jacobi operators.
\end{enumerate}
See \cite{GRKVL} for a classification of locally symmetric four-dimensional Osserman manifolds and \cite{Calvaruso-Zaeim14} for a description of gradient Ricci solitons on non-reductive homogeneous spaces, where metrics of Type (A.1) with $b=0$ play a distinguished role. 
\end{remark}

\section{Conformally Einstein non-reductive homogeneous spaces. The proof of the Main Theorem}\label{se:4}

The purpose of this section is to prove the Main Theorem, determining which non-reductive homogeneous four-manifolds contain an Einstein metric in their conformal class. We will exclude from our analysis the trivial cases of Einstein and locally conformally flat manifolds. Moreover, we will obtain the explicit form of the conformal Einstein metric.
Since any conformally Einstein manifold is necessarily Bach flat, Theorem \ref{th: Bach flat} shows that the analysis of the conformally Einstein equation  
\begin{equation}
\label{3--1}
2\Hes_\varphi +\varphi\,\rho= \frac{1}{4}\{2\Delta\varphi+\varphi\,\tau\}g
\end{equation}
must be carried out only for the following cases:
\begin{enumerate}
\item Type $(A.1)$ with $q=0$ and $b\neq 0$, or $q=-\frac{3a}{4}$ and $b\neq 0$.
\item Type $(A.2)$ with $\alpha=1$ and $b\neq 0$.
\item Type $(A.3)$ with $\epsilon=\pm 1$ and $b\neq\mp q$.
\item Type $(B.1)$ with $q=0\neq c$.
\end{enumerate}

It is important to emphasize that although any locally conformally Einstein metric is Bach flat, there are examples of strictly Bach flat manifolds, i.e., they are neither half conformally flat nor locally conformally Einstein (see, for example \cite{AGS, nuevo, LNu} and references therein). Indeed, one has the following necessary conditions for any solution of \eqref{3--1} (see also \cite{GoNu}).

\begin{proposition}\cite{Kozameh-Newman-Tod85}\label{necessary condition}
Let $(M,g)$ be a four-dimensional pseudo-Riemannian manifold such that $\bar{g}=e^{2\sigma}g$ is Einstein. Then 
\begin{enumerate}
\item $\mathfrak{C}+W(\cdot,\cdot,\cdot,\nabla \sigma)=0$,
\item $\bach=0$,
\end{enumerate}
where $\mathfrak{C}$ and $\bach$ are the Cotton and the Bach tensor, respectively.
\end{proposition}

Recall that the solutions $\varphi$ of the conformally Einstein equation \eqref{3--1} and the functions $\sigma$ in Proposition \ref{necessary condition}-(1) are related by $\sigma=-2 \log(\varphi)$.  Also, as a matter of notation, define a $(0,3)$-tensor field $\mathcal{C}$ by $\mathcal{C}=\mathfrak{C}+W(\cdot,\cdot,\cdot,\nabla \sigma)$. Obviously, $\mathcal{C}_{ijk}=-\mathcal{C}_{jik}$ for all $i,j,k\in\{1,\dots,4\}$.

Conditions $(1)-(2)$ above are also sufficient to be conformally Einstein if $(M,g)$ is \emph{weakly-generic} (i.e., the Weyl tensor, viewed as a map $TM\rightarrow\bigotimes^3 TM$ is injective). Note that cases $(1)-(3)$ in Theorem \ref{th: Bach flat} are not weakly-generic and thus we must study the existence of solutions of Equation \eqref{3--1} case by case. In opposition, metrics corresponding to Theorem \ref{th: Bach flat}-$(4)$ are weakly-generic.

\subsection{Type (A.1) with $q=0$ and $b\neq 0$, or $q=-\frac{3a}{4}$ and $b\neq 0$.}
We consider the two possibilities separately.
\subsubsection{Type (A.1) with $q=0$ and $b\neq 0$}
In this case by Equation \eqref{eq: Cotton A.1} the non-zero components of the Cotton tensor are given by
\begin{equation}
\mathfrak{C}_{121}=-\frac{24\,  b x_2}{a}\,,\qquad \mathfrak{C}_{122}=-\frac{12\,  b}{a}\,,
\end{equation}
and, by Equation \eqref{eq: Weyl A.1}, the only non-zero component of the Weyl tensor is given by
\begin{equation}
W_{1212}=-6 b\,,
\end{equation} 
which shows that $(M,g)$ is not weakly-generic.
For an arbitrary positive function $\varphi(x_1,x_2,x_3,x_4)$ on $M$, let $\sigma=-2 \log(\varphi)$.
Then a straightforward calculation shows that the gradient of $\sigma$ is given in the coordinate basis by 
\[
\begin{array}{l}
\nabla \sigma=\frac{4}{a^2\varphi} \left\{  \left(a x_4-c\right) \varphi_4+a   \varphi_3\right\}\partial_1\\
\noalign{\medskip}
\phantom{\nabla \sigma=}
-\frac{2}{a^2\varphi} \left\{ \varphi_4 \left(a+4 x_2 \left(a x_4-c\right)\right)+4 a x_2\varphi_3\right\}\partial_2\\
\noalign{\medskip}
\phantom{\nabla \sigma=}
+\frac{4}{a^2\varphi} \left\{  a \left(2 \varphi_3-2 x_2
   \varphi_2+\varphi_1\right)+2 \varphi_4 \left(a
   x_4-c\right)\right\}\partial_3\\
\noalign{\medskip}
\phantom{\nabla \sigma=}
+\frac{2}{a^3\varphi}\left\{ \varphi_4 \left(a b+4 a x_4 \left(a x_4-2 c\right)+4 c^2\right)+2a\varphi_1 \left(a x_4-c\right)\right.\\
\noalign{\medskip}
\phantom{\nabla \sigma=+\frac{2}{a^3\varphi}}
+\left.a \left(4\varphi_3 \left(a x_4-c\right)+\varphi_2 \left(4 x_2 \left(c-a x_4\right)-a\right)
\right)\right\}\partial_4\,,
\end{array}
\]
where $\partial_i=\frac{\partial}{\partial {x_i}}$ are the coordinate vector fields and $\varphi_i=\frac{\partial}{\partial {x_i}}\varphi$ denote the corresponding partial derivatives. 

Thus, the only non-zero components of the tensor $\mathcal{C}=\mathfrak{C}+W(\cdot,\cdot,\cdot,\nabla \sigma)$ are those given by
\begin{equation}\label{eq:A.1 q0}
\begin{array}{l}
a^2\varphi\mathcal{C}_{121}=-12 b \left(\varphi_4 \left(a-4 x_2 \left(c-a x_4\right)\right)+4 a x_2 \varphi_3+2 a x_2
   \varphi\right)\,,\\
\noalign{\medskip}
a^2\varphi\mathcal{C}_{122}=-12 b \left(-2 \varphi_4 \left(c-a x_4\right)+2 a \varphi_3+a
   \varphi\right)\,.
\end{array}
\end{equation}
Since $\mathcal{C}=0$ is a necessary condition for $(M,g)$ to be conformally Einstein,
 $a\varphi(\mathcal{C}_{121}-2x_2\mathcal{C}_{122})=-12 b\varphi_4$
 must be zero and, since $b\neq 0$,  in this case $\varphi$ does not depend on the coordinate $x_4$. Then
\[
\mathcal{C}_{122}=\frac{-12b(\varphi+2\varphi_3)}{a\varphi}\,,
\quad
\mathcal{C}_{121}=2x_2\mathcal{C}_{122}\,.
\]
Hence, $\mathcal{C}=0$ shows that  
\begin{equation}
\label{funcion-1-1}
\varphi(x_1,x_2,x_3)=e^{-\frac{x_3}{2}}\phi(x_1,x_2)\,,
\end{equation}
for some smooth function $\phi(x_1,x_2)$.

Now, we analyze the existence of solutions of \eqref{3--1} for some $\varphi$ as above. In order to simplify the notation, set 
\[
\mathcal{E}=2 \Hes_\varphi +\varphi\,\rho- \frac{1}{4}\{2\Delta\varphi+\varphi\,\tau\}g\,,
\]
and determine the conditions for $\mathcal{E}=0$.

Since $\mathcal{E}(\partial_1,\partial_1)=2e^{-\frac{x_3}{2}}\phi_{11}(x_1,x_2)$, any solution of \eqref{3--1} must be of the form \eqref{funcion-1-1} with $\phi(x_1,x_2)=\alpha_1(x_2)+x_1\alpha_2(x_2)$ for some smooth  functions $\alpha_1$, $\alpha_2$ on $M$. A calculation of  $\mathcal{E}(\partial_1,\partial_2)=-2 e^{-\frac{x_3}{2}} (\alpha_1'(x_2) + 
(x_1-1) \alpha_2'(x_2))$, shows that  $\alpha_1(x_2)=\kappa_1$ and $\alpha_2(x_2)=\kappa_2$ for some constants $\kappa_1$, $\kappa_2$. Further, the component $\mathcal{E}(\partial_2,\partial_4)=-2 \kappa_2 e^{-\frac{x_3}{2}}$ shows that $\kappa_2=0$ and hence \eqref{funcion-1-1} reduces to
\[
\varphi=\kappa_1 e^{-\frac{x_3}{2}}\,.
\] 
Now, a straightforward calculation shows that $\mathcal{E}=0$ holds and the conformal metric $\bar g=\varphi^{-2}g$ is Ricci flat. 

\begin{remark}
\rm
Since any non-reductive homogeneous manifold of Type $(A.1)$ with $q=0$ and $b\neq 0$ is conformally Osserman with two-step nilpotent conformal Jacobi operators, and this property is conformally invariant, the metric  $\bar g$ is Osserman with two-step nilpotent Jacobi operators.
\end{remark}

\subsubsection{Type (A.1) with $q=-\frac{3a}{4}$ and $b\neq 0$.}
We proceed as in the previous case.
For an arbitrary positive function $\varphi(x_1,x_2,x_3,x_4)$ on $M$, consider  $\sigma=-2 \log(\varphi)$.
Then  
\[
\begin{array}{l}
\nabla \sigma
=\frac{1}{a^2\varphi}\left\{ a \left(\varphi_3+3x_2 \varphi_2-\frac{3}{2}
   \varphi_1\right)+\varphi_4 \left(a x_4-c\right)\right\}\partial_1\\
\noalign{\medskip}
\phantom{\nabla \sigma =}
-\frac{1}{a^2\varphi}\left\{2 \varphi_4 \left(x_2 \left(a x_4-c\right)+a\right)+a x_2 \left(2
   \varphi_3+6 x_2 \varphi_2-3
   \varphi_1\right)\right\}\partial_2\\
\noalign{\medskip}
\phantom{\nabla \sigma =}
+\frac{1}{a^2\varphi}\left\{ a \left(2 \varphi_3-2 x_2
   \varphi_2+\varphi_1\right)+2 \varphi_4
   \left(a x_4-c\right)\right\}\partial_3\\
\noalign{\medskip}
\phantom{\nabla \sigma =}
+\frac{1}{a^3\varphi}\left\{ 2 \varphi _4 \left(a^2 x_4^2+a b-2 a c x_4+c^2\right)\right.\\
\noalign{\medskip}
\phantom{...............}
\left.
+ a \left(\varphi _1 \left(a x_4-c\right)- 2 \left(\varphi _2 \left(x_2 \left(a x_4-c\right)+a\right)+\varphi _3 \left(c-a x_4\right)\right)\right)\right\}\partial_4\,.
\end{array}
\]
Recall from \eqref{eq: Cotton A.1} that  the non-zero components of the Cotton tensor are given by
\begin{equation*}
\begin{array}{l}
\mathfrak{C}_{121}=\frac{12 b x_2}{a}\,,\qquad \mathfrak{C}_{122}=\frac{6 b}{a}\,,\qquad \mathfrak{C}_{131}= -\frac{12 b x_2^2}{a}\,,\\
\noalign{\medskip}
\mathfrak{C}_{232}=-\frac{3 b}{a}\,,
\qquad
\mathfrak{C}_{132}=\mathfrak{C}_{231}=2x_2\,\mathfrak{C}_{232}\,.
\end{array}
\end{equation*}
Equation \eqref{eq: Weyl A.1} shows that the non-zero components of the Weyl tensor are
\[
\begin{array}{lll}
W_{1212}=-3 b\,,& W_{1213}=3 b x_2\,,& W_{1223}=\frac{3 b}{2}\,,\\
\noalign{\medskip}
W_{1313}=-3 b x_2^2,& W_{1323}=-\frac{3 b x_2}{2},& W_{2323}=-\frac{3 b}{4}\,.
\end{array}
\]
Then, the non-zero components of the tensor field $\mathcal{C}$ are given by
\[
\begin{array}{l}
\mathcal{C}_{131}=-x_2\mathcal{C}_{121},\quad \mathcal{C}_{231}=-\frac{1}{2}\mathcal{C}_{121}\,,\\
\noalign{\medskip}
\mathcal{C}_{132}=-x_2\mathcal{C}_{122},\quad \mathcal{C}_{232}=-\frac{1}{2}\mathcal{C}_{122}\,,\\
\noalign{\medskip}
\mathcal{C}_{133}=-x_2\mathcal{C}_{123},\quad \mathcal{C}_{233}=-\frac{1}{2}\mathcal{C}_{123}\,,
\end{array}
\]
where
\[
\begin{array}{l}
a^2\varphi\,\mathcal{C}_{121}=6 b \left(x_2 \left(2 a \varphi +a \varphi _1-2 a \varphi _3-2 a \varphi _4 x_4+2 c \varphi _4\right)-a \varphi _4-2 a \varphi _2 x_2^2\right),\\
\noalign{\medskip}
a^2\varphi\,\mathcal{C}_{122}=3 b \left(2 a \varphi +a \varphi _1-2 a \varphi _3-2 a \varphi _2 x_2-2 a \varphi _4 x_4+2 c \varphi _4\right)\,,\\
\noalign{\medskip}
a^2\varphi\,\mathcal{C}_{123}=-3 a b \varphi _4\,.
\end{array}
\]

Since $a,b\neq0$ and $\mathcal{C}_{123}=0$ the function $\varphi(x_1,x_2,x_3,x_4)$ does not depend on the coordinate $x_4$ and the tensor field $\mathcal{C}_{121}$ reduces to  
$$
\mathcal{C}_{122}=\frac{3b}{a\varphi}\left\{ 
2 \varphi + \varphi _1-2  \varphi _3-2  \varphi _2 x_2\right\}\,,
\quad
\mathcal{C}_{121}=2x_2 \mathcal{C}_{122}\,.
$$
A solution of the differential equation 
$2 \varphi = 2  \varphi _3+2 x_2 \varphi _2 -\varphi _1$
is necessarily of the form
\begin{equation}
\label{funcion-1-2}
\varphi(x_1,x_2,x_3)=e^{-2 x_1} \phi(e^{2x_1}x_2,2x_1+x_3)
=e^{-2 x_1} (\phi\circ\psi)(x_1,x_2,x_3)\,,
\end{equation}
where $\psi(x_1,x_2,x_3)=(e^{2x_1}x_2,2x_1+x_3)$ and $\phi(z,\omega)$ is an arbitrary function
for $z=e^{2x_1}x_2$ and $\omega=2x_1+x_3$.

Now, we analyze the existence of solutions of \eqref{3--1} for some $\varphi$ as in \eqref{funcion-1-2}. Setting 
\[
\mathcal{E}=2 \Hes_\varphi +\varphi\,\rho- \frac{1}{4}\{2\Delta\varphi+\varphi\,\tau\}g\,,
\]
one has $\mathcal{E}(\partial_2,\partial_2)=2e^{2x_1}\partial^2_{z^2}\phi=0$, and hence
$$
\varphi(x_1,x_2,x_3)=e^{-2 x_1} \left( e^{2x_1}x_2\hat{\phi}(2x_1+x_3)+\bar{\phi}(2x_1+x_3)  \right)
$$
for some smooth functions $\hat{\phi}(\omega)$, $\bar{\phi}(\omega)$. Considering now the component 
$\mathcal{E}(\partial_2,\partial_3)=2\hat{\phi}'-\hat{\phi}=0$, one has that
$\hat{\phi}(2x_1+x_3)=\kappa e^{\frac{1}{2}(2x_1+x_3)}$ for some constant $\kappa$. Now the only non-zero components of the tensor field $\mathcal{E}$ are given by
$$
\begin{array}{l}
\mathcal{E}(\partial_1,\partial_1)=2\mathcal{E}(\partial_1,\partial_3)=2\mathcal{E}(\partial_3,\partial_3)\\
\noalign{\medskip}
\phantom{\mathcal{E}(\partial_1,\partial_1)}
= 4e^{-2x_1}\left(\bar{\phi}-3\bar{\phi}'+2\bar{\phi}''   \right)\,.
\end{array}
$$
Hence $\mathcal{E}=0$ gives $\bar{\phi}(2x_1+x_3)=\kappa_1 e^{\frac{1}{2}(2x_1+x_3)}+\kappa_2 e^{2x_1+x_3}$ and thus any solution of the conformally Einstein equation is of the form
$$
\varphi(x_1,x_2,x_3,x_4)=\kappa_1 e^{x_3} + \kappa_2 e^{\frac{1}{2}x_3 - x_1} + \kappa_3\, x_2\, e^{\frac{1}{2}x_3 + x_1}\,.
$$
Moreover, any of the conformal metrics $\bar g=\varphi^{-2}g$ is Ricci flat.

\begin{remark}
\rm
Since any non-reductive homogeneous manifold of Type (A.1) with $q=-\frac{3}{4}a$ and $b\neq 0$ is conformally Osserman with two-step nilpotent conformal Jacobi operators, any conformal Einstein metric $\bar g$ is Osserman with two-step nilpotent Jacobi operators. Moreover, there is a $3$-parameter family of conformally equivalent Osserman metrics.  
This shows that the cases $q=0$ and $q=-\frac{3}{4}a$ are essentially different since the space of conformally Einstein metrics is one-dimensional in the first case and three-dimensional in the second one.
\end{remark}

\subsection{Type (A.2) with $\alpha=1$ and $b\neq 0$}
Let $\varphi(x_1,x_2,x_3,x_4)$ be a positive function on $M$ and $\sigma=-2 \log(\varphi)$.
Then 
\[
\begin{array}{l}
\nabla \sigma= 
\frac{2}{a^2 q \varphi}
\left\{ a e^{-2 x_4} \left(q \varphi _3-c \varphi _4\right)- e^{-4x_4} \varphi _1 \left(c^2-b q\right)\right\}\partial_1\\
\noalign{\medskip}
\phantom{\nabla \sigma=}
-\frac{2}{a \varphi} \varphi _2 e^{-2 x_4}\,\partial_2
+\frac{2}{a \varphi} \varphi _1 e^{-2 x_4}\,\partial_3
-\frac{1}{a q \varphi}\left\{ 2 a \varphi _4+2  c \varphi _1 e^{-2 x_4}\right\}\partial_4\,.
\end{array}
\]

It follows from \eqref{eq: Cotton A.2} and \eqref{eq: Weyl A.2} that the non-zero components of the Cotton and the Weyl tensors are given by
$$
\mathfrak{C}_{343}=\frac{b}{q}\,,  \qquad \text{and} \qquad
W_{2323}=\frac{a b }{2 q}e^{2 x_4}\,,\quad W_{3434}=-\frac{b}{2}\,,
$$
respectively, from where it follows that $(M,g)$ is not weakly-generic. Therefore,  the only non-zero components of the tensor field 
$\mathcal{C}=\mathfrak{C}+W(\cdot,\cdot,\cdot,\nabla \sigma)$ are those given by
\begin{equation}\label{eq:A.2}
\begin{array}{l}
\mathcal{C}_{232}=-\frac{b}{q}\frac{\varphi_1}{\varphi}\,,\qquad \mathcal{C}_{233}=-\frac{b}{q}\frac{\varphi_2}{\varphi}\,,\qquad
\mathcal{C}_{344}=-\frac{b}{a}\frac{\varphi_1}{\varphi}e^{-2 x_4}\,,
\\
\noalign{\medskip}
\mathcal{C}_{343}=\frac{b}{a \,q\, \varphi} \left(a \left(\varphi -\varphi _4\right)-c \varphi _1 e^{-2 x_4}\right)\,.
\end{array}
\end{equation}
Since, $a\, b\neq 0$  the first two equations in \eqref{eq:A.2} show that $\varphi(x_1,x_2,x_3,x_4)$ does not depend on the coordinates $x_1$ and $x_2$. Hence, the tensor field $\mathcal{C}$ reduces to
\[
\mathcal{C}_{343}=\frac{b \left(\varphi -\varphi _4\right)}{q \varphi }\,,
\]
where $\varphi$ is a smooth function on the coordinates $(x_3,x_4)$ and it follows from $\mathcal{C}_{343}=0$ that $\varphi(x_3,x_4)=\phi(x_3)e^{x_4}$, for some smooth function $\phi(x_3)$. 

Considering now the conformally Einstein equation \eqref{3--1}, and setting
\[
\mathcal{E}=2 \Hes_\varphi +\varphi\,\rho- \frac{1}{4}\{2\Delta\varphi+\varphi\,\tau\}g\,,
\]
the only non-zero component of the tensor $\mathcal{E}$ is 
\[
\mathcal{E}(\partial_3,\partial_3)=\frac{e^{x_4} \left(2 q \phi ''-b \phi \right)}{q}\,.
\]
Integrating $\mathcal{E}(\partial_3,\partial_3)=0$ we obtain that 
\begin{equation}\label{A2-final}
\left\{
\begin{array}{lr}
\varphi=e^{x_4-x_3\sqrt{\frac{b}{2 q}}}\left(\kappa_1e^{x_3\sqrt{\frac{2b}{q}}}+\kappa_2\right),& \text{if}\quad  b\,q>0,\\
\noalign{\medskip}
\varphi=e^{x_4}\left(\kappa_1 \cos\left(x_3\sqrt{-\frac{b}{2q}}\right)+\kappa_2 \sin\left(x_3\sqrt{-\frac{b}{2q}}\right)\right),& \text{if}\quad  b\,q<0\,.
\end{array}
\right.
\end{equation}
Moreover, a long but straightforward computation shows that the metric $\bar g=\varphi^{-2}g$ for any function $\varphi$ given by \eqref{A2-final} is Ricci flat.

\begin{remark}
\label{re:A2-conf-Einstein}\rm
For each of the possibilities in \eqref{A2-final} there are at least two conformal metrics which are Einstein (indeed, Ricci flat). Moreover, for any of the conformal Einstein metrics, there are some conformal deformation of the metric which remains Einsteinian. 

Further observe that no metric (A.2) with $\alpha=1$ and $b\neq 0$ is half conformally flat, and hence they are not in the conformal class of any Osserman metric.
\end{remark}

\subsection{Type (A.3) with $\epsilon=\pm 1$ and $b\neq\mp q$}

We will briefly schedule the proof of the case corresponding to $\epsilon=1$. The situation when $\epsilon=-1$ is completely analogous. Hence assume $b\neq-q$. As in the previous cases, let $\varphi(x_1,x_2,x_3,x_4)$ be a positive function and set $\sigma=-2 \log(\varphi)$. Then
\[
\begin{array}{l}
\nabla \sigma=\frac{1}{a^2 \, b \varphi }\left\{2 e^{-4 x_3} \left(a e^{2 x_3} \left(c \varphi _3-b \varphi _4\right)+\varphi _1 \left(b q-c^2\right)\right)\right\}\partial_1
\\
\noalign{\medskip}
\phantom{=}
-\frac{2\varphi _2}{a\varphi }\left\{e^{-2 x_3} \sec\left(x_4\right)^2\right\}\partial_2
+\frac{2}{a\, b \varphi } \left\{ c \varphi _1 e^{-2 x_3}-a \varphi _3\right\}\partial_3
-\frac{2}{a\varphi } \varphi _1 e^{-2 x_3}\partial_4\,.
\end{array}
\]
It follows  from  \eqref{eq: Cotton A.31} and \eqref{eq: Weyl A.31} that the non-zero components of the tensor $\mathcal{C}=\mathfrak{C}+W(\cdot,\cdot,\cdot,\nabla \sigma)$ are given by
\begin{equation}\label{eq:A.31}
\begin{array}{l}
b\, \varphi\,\mathcal{C}_{242}=  (b+q) \cos\left(x_4\right)^2\varphi _1\,, \quad b\, \varphi\,\mathcal{C}_{244}=- (b+q)\varphi _2 \,,\\
\noalign{\medskip}
a\varphi\,\mathcal{C}_{343}= -(b+q)e^{-2 x_3} \varphi _1\,,\\
\noalign{\medskip}
a\,b\,\varphi\,\mathcal{C}_{344}=-(b+q)e^{-2 x_3}  \left(a \left(\varphi -\varphi _3\right) e^{2 x_3}+c \varphi _1\right)\,.
\end{array}
\end{equation}
Since $a\, \neq 0$ and $b\neq -q$ the first two equations in \eqref{eq:A.31} show that $\varphi$ does not depend on the coordinates $x_1$ and $x_2$ and the tensor field $\mathcal{C}$ reduces to 
\[
b\,\varphi\,\mathcal{C}_{344}=-(b+q)\left(\varphi -\varphi _3\right)\,,
\]
where $\varphi$ is a smooth function on the coordinates $(x_3,x_4)$. Now $\mathcal{C}_{344}=0$ gives $\varphi(x_3,x_4)=\phi(x_4)e^{x_3}$, for some smooth function $\phi(x_4)$. 

Consider now the conformally Einstein equation and set, as in the previous cases,
$
\mathcal{E}=2 \Hes_\varphi +\varphi\,\rho- \frac{1}{4}\{2\Delta\varphi+\varphi\,\tau\}g\,.
$
A straightforward calculation shows that the only non-zero component of the tensor field $\mathcal{E}$ is given by
\[
\mathcal{E}(\partial_4,\partial_4)=\frac{1}{b}e^{x_3} \left((b-q)\phi  +2 b \phi ''\right)\,,
\]
which shows that  $\phi(x_4)$ is determined by the equation
$\phi ''=-\frac{b-q}{2b}\phi$. Hence the conformal deformation $\varphi(x_3,x_4)$ is given by 
\begin{equation}\label{final +}
\left\{ 
\begin{array}{ll}
\!\!\!\varphi_+=
(\kappa_1 x_4+\kappa_2)e^{x_3},&  \text{if}\,\, b-q=0,\\
\noalign{\medskip}
\!\!\!\varphi_+=
e^{x_3-x_4\sqrt{\frac{q-b}{2 b}}}\left(\kappa_1e^{x_4\sqrt{\frac{2(q-b)}{b}}}  +\kappa_2\right),& \text{if}\,\, b(q-b)>0,\\
\noalign{\medskip}
\!\!\!\varphi_+=
e^{x_3}\!\left(\!\kappa_1\cos\!\left(\!x_4\sqrt{\frac{b-q}{2b}}\right)\!+\kappa_2\sin\!\left(\!x_4\sqrt{\frac{b-q}{2b}}\right)\!\!\right),
& \text{if}\,\, b(q-b)<0\,.
\end{array}
\right.
\end{equation}
Moreover in all the cases above the conformal metric $\bar g_+=\varphi_+^{-2}g_+$ is Ricci flat.

The case $\epsilon=-1$ is obtained in a completely analogous way. For any metric $g_-$ given by \eqref{eq:A.3-2}, the conformal metric $\bar g_-=\varphi_-^{-2}g_-$ is Ricci flat, where 

\begin{equation}\label{final -}
\left\{
\begin{array}{ll}
\!\!\!\varphi_-\!=\!(\kappa_1 x_4+\kappa_2)e^{x_3},&\! \text{if}\,\, b+q=0,\\
\noalign{\medskip}
\!\!\!\varphi_-\!=\!e^{x_3-x_4\sqrt{\frac{q+b}{2 b}}}\left(\kappa_1e^{x_4\sqrt{\frac{2(q+b)}{b}}}\!  +\kappa_2\right),& \!\text{if}\,\,  b(q+b)>0,\\
\noalign{\medskip}
\!\!\!\varphi_-\!=\!e^{x_3}\!\!\left(\!\kappa_1\cos\!\left(\!x_4\sqrt{\!-\frac{b+q}{2b}}\right)\!+\kappa_2\sin\!\left(\!x_4\sqrt{\!-\frac{b+q}{2b}}\right)\!\!\right),
 &\!  \text{if}\,\, b(q+b)<0\,.
\end{array}
\right.
\end{equation}

\begin{remark}
\label{re:A3-conf-Einstein}\rm
For each of the possibilities in \eqref{final +} and \eqref{final -}  there are at least two conformal metrics which are Einstein. Equivalently, for any conformally Einstein metric, there are some conformal deformation of the metric which remains to be Einstein. 

Further observe that no metric (A.3) with $\epsilon=\pm 1$ and $b\neq\mp q$ is half conformally flat, and hence $(M,g)$ is not in the conformal class of an Osserman metric.
\end{remark}

\subsection{Type (B.1) with $q=0\neq c$}
Setting $q=0$ in Equations \eqref{eq: Cotton B.1} and \eqref{eq: Weyl B.1}, the non-zero components of the Cotton and the Weyl tensors are given by
$$
\begin{array}{lll}\mathfrak{C}_{121}=\frac{45 c^2 x_2}{a^2}\,,& \mathfrak{C}_{122}=\frac{45 c^2}{2 a^2}\,, &\text{and}
\\
\noalign{\bigskip}
W_{1212}=\frac{3 c^2 x_3}{a}-3 \left(b+2 c x_4\right)\,,& 
W_{1213}=-\frac{3 c \left(a+2 c x_2\right)}{2 a}\,,&\\
\noalign{\medskip}
W_{1214}=-3 c x_2\,,&
W_{1223}=-\frac{3 c^2}{2 a}\,,& 
W_{1224}=-\frac{3 c}{2}\,,
\end{array}
$$
respectively. This shows that, in opposition to the previous cases, $(M,g)$ is weakly-generic and thus 
$\mathcal{C}=\mathfrak{C}+W(\cdot,\cdot,\cdot,\nabla \sigma)=0$ is a necessary and sufficient condition to be conformally Einstein.

As in the previous cases, consider  $\varphi(x_1,x_2,x_3,x_4)$  a positive function and set $\sigma=-2\log(\varphi)$. Express the gradient of $\sigma$ as
\[
\begin{array}{l}
\nabla \sigma= \frac{2}{a^2\,\varphi}\left\{c\varphi_4-a\varphi_3\right\}\partial_1
+\frac{2}{a^2\,\varphi}\left\{2 x_2 \left(a \varphi _3-c \varphi _4\right)- a \varphi _4\right\}\partial_2\\
\noalign{\medskip}
\phantom{\nabla \sigma=}
+\frac{2}{a^2\,\varphi}\left\{ x_3 \left(2 a \varphi _3-c \varphi _4\right)- a \varphi _1+2 a \varphi _2 x_2\right\}\partial_3\\
\noalign{\medskip}
\phantom{\nabla \sigma=}

+\frac{2}{a^2\,\varphi}\left\{b \varphi _4-\varphi _2 \left(a+2 c x_2\right)+c \varphi _1-c \varphi _3 x_3+2 c \varphi _4 x_4\right\}\partial_4\,.
\end{array}
\]

Now, the components $\mathcal{C}_{123}$ and $\mathcal{C}_{124}$ of the tensor field $\mathcal{C}=\mathfrak{C}+W(\cdot,\cdot,\cdot,\nabla \sigma)$ are given by
$$
a\,\varphi\,\mathcal{C}_{123}=3 c \varphi _3,\quad a\,\varphi\,\mathcal{C}_{124}=3 c \varphi _4\,,
$$
and, since $c\, \neq 0$ and $\mathcal{C}_{123}=\mathcal{C}_{124}=0$, the function $\varphi$ is independent of the coordinates $x_3$ and $x_4$. Assuming $\varphi$ to be a smooth function on the coordinates $(x_1,x_2)$, the non-zero components of $\mathcal{C}$ reduce to
\begin{equation}
\mathcal{C}_{121}=-\frac{3 c \left(a \varphi _1-15 c \varphi  x_2\right)}{a^2 \varphi }\,,\quad \mathcal{C}_{122}=\frac{3 c \left(15 c \varphi -2 a \varphi _2\right)}{2 a^2 \varphi }\,.
\end{equation}
A straightforward computation shows that $\mathcal{C}_{122}=0$ if and only if 
\begin{equation}
\varphi=\phi(x_1) e^{\frac{15 c}{2 a} x_2}\,,
\end{equation}
for some smooth function $\phi(x_1)$. Then $\mathcal{C}_{121}$ becomes
\begin{equation}
\mathcal{C}_{121}=-\frac{3 c \left(a \phi'(x_1)-15 c x_2 \phi(x_1) \right)}{a^2 \phi(x_1) }\,,
\end{equation}
from where it follows that $\phi$ vanishes identically, and hence $\varphi\equiv 0$, which is a contradiction. Hence this manifold is not conformally Einstein.

\begin{remark}
\rm
Observe that the conformally Einstein metrics in the Main Theorem--(1) are always of neutral signature, while metrics corresponding to cases (2) and (3) may be either Lorentzian or of neutral signature $(2,2)$, depending on the choice of the parameters defining the metrics \eqref{eq:A.2-1}, \eqref{eq:A.3-1} and \eqref{eq:A.3-2}.
\end{remark}


\begin{thebibliography}{99}
\bibitem{AGS}
E. Abbena, S. Garbiero, and S. Salamon, 
Bach-flat Lie groups in dimension $4$,
\emph{C. R. Math. Acad. Sci. Paris} \textbf{351} (2013), 303--306. 

\bibitem{BGi}
N. Bla\v{z}i\'c and P. Gilkey, 
Conformally Osserman manifolds and self-duality in Riemannian geometry, 
\emph{Differential geometry and its applications}, 15–-18, Matfyzpress, Prague, 2005. 

\bibitem{Brinkmann24}
H.W. Brinkmann, Riemann spaces conformal to Einstein spaces, \emph{Math. Ann.} \textbf{91} (1924), 269--278.

\bibitem{Brinkmann25}
H.W. Brinkmann, Einstein spaces which are mapped conformally on each other, \emph{Math. Ann.}  \textbf{94} (1925),  119--145.

\bibitem{Calvaruso-Fino}
G. Calvaruso and A. Fino, Ricci Solitons and Geometry of Four-dimensional Non-reductive
Homogeneous Spaces, \emph{Canad. J. Math.} \textbf{64} (2012), 778--804.

\bibitem{Calvaruso-Fino-Zaeim15}
G. Calvaruso, A. Fino, and A. Zaeim, Homogeneous geodesics of non-reductive homogeneous pseudo-Riemannian $4$-manifolds, \emph{Bull. Brazil. Math. Soc.} \textbf{46} (2015), 1--42.

\bibitem{Calvaruso-Zaeim14}
G. Calvaruso and A. Zaeim, A complete classification of Ricci and Yamabe solitons of non-reductive homogeneous $4$-spaces, \emph{J. Geom. Phys.} \textbf{80} (2014), 15–-25.

\bibitem{Calvaruso-Zaeim14-2}
G. Calvaruso and A. Zaeim, Geometric structures over non-reductive homogeneous 4-spaces, Adv. Geom. 14 (2014), 191--214.


\bibitem{CLGRGVL}
E. Calvi\~no-Louzao, E. Garc\'\i a-R\'\i o, P. Gilkey, and R. V\'azquez-Lorenzo,
The geometry of modified Riemannian extensions, 
\emph{Proc. R. Soc. Lond. Ser. A Math. Phys. Eng. Sci.} \textbf{465} (2009), 2023--2040. 

\bibitem{nuevo}
E. Calvi\~no-Louzao, I. Guti\'errez-Rodr\'\i guez, E. Garc\'\i a-R\'\i o, and R. V\'azquez-Lorenzo,
New examples of strictly Bach-flat four-manifolds, 
to appear.

\bibitem{Fels-Renner06}
M. E. Fels and A. G. Renner, 
Non-reductive Homogeneous Pseudo-Riemannian Manifolds of Dimension Four, 
\emph{Canad. J. Math.} \textbf{58} (2006), 282–-311.

\bibitem{GaOu}
P. M. Gadea and J. A. Oubi\~na, 
Reductive homogeneous pseudo-Riemannian manifolds,
\emph{Monatsh. Math.} \textbf{124} (1997), 17--34.

\bibitem{GRKVL}
E. Garc\'\i a--R\'\i o, D. N. Kupeli, and R. V\'azquez-Lorenzo,
\emph{Osserman manifolds in semi-Riemannian geometry}, Lect. Notes
Math. \textbf{1777}, Springer-Verlag, Berlin, Heidelberg, New York,
2002.

\bibitem{GoNu}
R. Gover and P. Nurowski,
Obstructions to conformally Einstein metrics in $n$ dimensions, 
\emph{J. Geom. Phys.} \textbf{56} (2006), 450--484.

\bibitem{Kozameh-Newman-Tod85}
C. N. Kozameh, E. T. Newman, and K. P. Tod, 
Conformal Einstein spaces. \emph{Gen. Rel. Grav.} \textbf{17} (1985), 343--352.

\bibitem{Kuhnel-Rademacher} 
W. K\"{u}hnel and H.B. Rademacher, 
Conformal transformations of pseudo-Riemannian manifolds,
\emph{Recent developments in pseudo-Riemannian geometry}, 261--298, ESI Lect. Math. Phys., Eur. Math. Soc., Zürich, 2008. 

\bibitem{KR2}
W. K\"{u}hnel and H.B. Rademacher, 
Conformally Einstein spaces revisited,
\emph{Pure and Appl. Differential Geom. PADGE 2012 - In Memory of Franki Dillen}, 161--167, Shaker Verlag, Aachen, 2013. 

\bibitem{LNu}
Th. Leistner and P. Nurowski,
Ambient metrics for n-dimensional pp-waves,
\emph{Comm. Math. Phys.} \textbf{296} (2010), 881--898.

\bibitem{Ni}
Y. Nikolayevsky, 
Conformally Osserman manifolds,
\emph{Pacific J. Math.} \textbf{245} (2010), 315--358.
    
\bibitem{Walker}
A. G. Walker, Canonical form for a Riemannian space with
a parallel field of null planes, {\it Quart. J. Math. Oxford} {\bf
(2) 1} (1950), 69--79.
\end{thebibliography}
\end{document}